\documentclass[11pt]{article}

\usepackage{amssymb,latexsym,amsmath}

\usepackage{graphicx}

\begin{document}

\newcommand{\bfi}{\bfseries\itshape}

\makeatletter

\@addtoreset{figure}{section}

\def\thefigure{\thesection.\@arabic\c@figure}

\def\fps@figure{h, t}

\@addtoreset{table}{bsection}

\def\thetable{\thesection.\@arabic\c@table}

\def\fps@table{h, t}

\@addtoreset{equation}{section}

\def\theequation{\thesubsection.\arabic{equation}}

\makeatother

\newtheorem{thm}{Theorem}[section]

\newtheorem{prop}[thm]{Proposition}

\newtheorem{lema}[thm]{Lemma}

\newtheorem{cor}[thm]{Corollary}

\newtheorem{defi}[thm]{Definition}

\newtheorem{rk}[thm]{Remark}

\newtheorem{exempl}{Example}[section]

\newenvironment{exemplu}{\begin{exempl}  \em}{\hfill $\surd$

\end{exempl}}

\newcommand{\comment}[1]{\par\noindent{\raggedright\texttt{#1}

\par\marginpar{\textsc{Comment}}}}

\newcommand{\todo}[1]{\vspace{5 mm}\par \noindent \marginpar{\textsc{ToDo}}\framebox{\begin{minipage}[c]{0.95 \textwidth}

\tt #1 \end{minipage}}\vspace{5 mm}\par}

\newcommand{\ea}{\mbox{{\bf a}}}

\newcommand{\eu}{\mbox{{\bf u}}}

\newcommand{\ueu}{\underline{\eu}}

\newcommand{\ueo}{\overline{u}}

\newcommand{\oeu}{\overline{\eu}}

\newcommand{\ew}{\mbox{{\bf w}}}

\newcommand{\ef}{\mbox{{\bf f}}}

\newcommand{\eF}{\mbox{{\bf F}}}

\newcommand{\eC}{\mbox{{\bf C}}}

\newcommand{\en}{\mbox{{\bf n}}}

\newcommand{\eT}{\mbox{{\bf T}}}

\newcommand{\eL}{\mbox{{\bf L}}}

\newcommand{\eR}{\mbox{{\bf R}}}

\newcommand{\eV}{\mbox{{\bf V}}}

\newcommand{\eU}{\mbox{{\bf U}}}

\newcommand{\ev}{\mbox{{\bf v}}}

\newcommand{\eve}{\mbox{{\bf e}}}

\newcommand{\uev}{\underline{\ev}}

\newcommand{\eY}{\mbox{{\bf Y}}}

\newcommand{\eK}{\mbox{{\bf K}}}

\newcommand{\eP}{\mbox{{\bf P}}}

\newcommand{\eS}{\mbox{{\bf S}}}

\newcommand{\eJ}{\mbox{{\bf J}}}

\newcommand{\eB}{\mbox{{\bf B}}}

\newcommand{\eH}{\mbox{{\bf H}}}

\newcommand{\leb}{\mathcal{ L}^{n}}

\newcommand{\eI}{\mathcal{ I}}

\newcommand{\eE}{\mathcal{ E}}

\newcommand{\hen}{\mathcal{H}^{n-1}}

\newcommand{\eBV}{\mbox{{\bf BV}}}

\newcommand{\eA}{\mbox{{\bf A}}}

\newcommand{\eSBV}{\mbox{{\bf SBV}}}

\newcommand{\eBD}{\mbox{{\bf BD}}}

\newcommand{\eSBD}{\mbox{{\bf SBD}}}

\newcommand{\ecs}{\mbox{{\bf X}}}

\newcommand{\eg}{\mbox{{\bf g}}}

\newcommand{\paromega}{\partial \Omega}

\newcommand{\gau}{\Gamma_{u}}

\newcommand{\gaf}{\Gamma_{f}}

\newcommand{\sig}{{\bf \sigma}}

\newcommand{\gac}{\Gamma_{\mbox{{\bf c}}}}

\newcommand{\deu}{\dot{\eu}}

\newcommand{\dueu}{\underline{\deu}}

\newcommand{\dev}{\dot{\ev}}

\newcommand{\duev}{\underline{\dev}}

\newcommand{\weak}{\stackrel{w}{\approx}}

\newcommand{\mild}{\stackrel{m}{\approx}}

\newcommand{\strong}{\stackrel{s}{\approx}}

\newcommand{\weakdown}{\rightharpoondown}

\newcommand{\opg}{\stackrel{\mathfrak{g}}{\cdot}}

\newcommand{\opunu}{\stackrel{1}{\cdot}}
\newcommand{\opdoi}{\stackrel{2}{\cdot}}

\newcommand{\opn}{\stackrel{\mathfrak{n}}{\cdot}}

\newcommand{\tr}{\ \mbox{tr}}

\newcommand{\Ad}{\ \mbox{Ad}}

\newcommand{\ad}{\ \mbox{ad}}

\renewcommand{\contentsname}{ }

\title{Existence and construction  of bipotentials for graphs of  multivalued 
laws}

\author{Marius Buliga\footnote{"Simion Stoilow" Institute of Mathematics of the 
Romanian Academy, PO BOX 1-764,014700 Bucharest, Romania, e-mail: 
Marius.Buliga@imar.ro } ,
 G\'ery de Saxc\'e\footnote{Laboratoire de 
  M\'ecanique de Lille, UMR CNRS 8107, Universit\'e des Sciences et 
Technologies de Lille,
 Cit\'e Scientifique, F-59655 Villeneuve d'Ascq cedex, France, 
e-mail: gery.desaxce@univ-lille1.fr} , Claude  Vall\'ee\footnote{Laboratoire de 
M\'ecanique des Solides, UMR 6610, UFR SFA-SP2MI, bd M. et P. Curie, 
t\'el\'eport 2, BP 30179, 86962 Futuroscope-Chasseneuil cedex, 
France, e-mail: vallee@lms.univ-poitiers.fr} }

\date{This version: 08.01.2007}

\maketitle

{\bf MSC-class:} 49J53; 49J52; 26B25

\begin{abstract}
Based on an extension of Fenchel inequality, bipotentials are  non smooth 
mechanics tools,  used to model  various non associative multivalued 
constitutive laws of dissipative materials (friction contact, soils, cyclic
plasticity of metals, damage).

Let $X$, $Y$ be dual locally convex spaces, with duality product $\langle \cdot , \cdot \rangle: X \times Y \rightarrow \mathbb{R}$.
 Given the graph $\displaystyle M \subset X\times Y$ of a multivalued law 
$\displaystyle T:X\rightarrow 2^{Y}$, 
we state a simple necessary and sufficient condition for the existence of a 
bipotential $b$ for which $M$ is the set of $(x, y)$ such that
$b(x, y) = \langle x, y\rangle$.  

If this condition  is fulfilled, we use convex lagrangian covers in order  to 
construct such a bipotential, generalizing a theorem due to Rockafellar, 
which states that a multivalued constitutive law admits a superpotential if and 
only if its graph is cyclically monotone. 
\end{abstract}

\maketitle

\section{Introduction}

The basic tools of the mechanics of continua are the kinematical compatibility
 and equilibrium local equations but they are not sufficient to describe the 
deformation and motion of the continuous media. Additional information must be 
given through the constitutive laws traducing the material behaviour. In its 
simplest form, a constitutive law is given by a graph collecting couples of dual
 variables resulting from experimental testing.

For many   physically relevant situations, 
the constitutive laws are multivalued, but also associated. The graph of the 
constitutive law is included in the graph of the subdifferential of a convex 
(and lower semi continuous) superpotential $\phi$. The constitutive law takes the 
form of a  differential inclusion, $ y \in \partial \phi (x)$. 
Any superpotential 
$\phi$ has a polar function $\phi^* $  satisfying a fundamental relation, 
Fenchel's inequality, $ \forall x,y, \: \phi(x)+\phi^{*}(y) \geq \langle x, 
y \rangle $. The constitutive law may be also written as 
$ x \in \partial \phi^{*} (y)$. 
In the literature,  this kind of materials are often called 
standard materials or generalized standard materials \cite{Halp JM 75}.

From the viewpoint of applications it is important to know  whether there 
exists a superpotential for a given non smooth graph and how to construct it. 
The answer to this question is provided by a famous theorem due to Rockafellar 
\cite{rocka} that 
ensures a graph admits a superpotential if and only if it is maximal cyclically 
monotone.

However, some of the constitutive laws are non-associated. 
They cannot be cast in the mould of the standard materials. To skirt this 
pitfall, a possible response, proposed first in \cite{saxfeng}, consists in 
constructing a function $b$ of two variables, bi-convex and satisfying an 
inequality generalizing Fenchel's one, $ \forall x,y, \: \ b (x,y) \geq 
\langle x, y \rangle $. We call it a bipotential. Physically, it represents 
the dissipation. In the case of associated constitutive laws the bipotential 
has the expression $b(x,y) = \phi(x)+\phi^{*}(y)$.

 As for the non associated constitutive laws which can be expressed with the 
help of bipotentials, they have the  form of an implicit relation between dual 
variables,  $ y \in \partial b(\cdot , y)(x) $. In Mechanics they are 
called  implicit, or weak, normality rules.  
The applications of bipotentials to Solid Mechanics are various: 
Coulomb's friction law \cite{sax CRAS 92} , non-associated Dr\"ucker-Prager 
\cite{sax boush KIELCE 93}  and Cam-Clay models \cite{sax BOSTON 95}  in 
Soil Mechanics, cyclic Plasticity (\cite{sax CRAS 92},\cite{bodo sax EJM 01}) 
and Viscoplasticity \cite{hjiaj bodo CRAS 00} of metals with non linear 
kinematical hardening rule, Lemaitre's damage law \cite{bodo}, the coaxial 
laws (\cite{dangsax},\cite{vall leri CONST 05}). Such kind of materials are 
called implicit standard materials. A synthetic review of these laws expressed 
in terms of bipotentials  can be 
found in \cite{dangsax} and \cite{vall leri CONST 05}.

The use of bipotentials in applications  is particularly attractive in 
numerical simulations when using the finite element method, but the interest is 
not limited to this aspects. For instance, the bound theorems of the limit 
analysis (\cite{sax boussh IJMS 98}, \cite{boush chaa IJMS 02}) and the plastic 
shakedown theory (\cite{sax trit AACHEN 00},  \cite{dangsax}, 
\cite{boush chaa IJP 03}, \cite{bouby sax IJSS 06}) can be reformulated in a 
broader framework, precisely by means of weak normality rules. 
From an applied numerical viewpoint, the bipotential method suggests new 
algorithms, fast but robust, as well as variational error estimators assessing 
the accurateness of the finite element mesh (\cite{hjiaj sax BUDA 96}, 
\cite{hjiaj sax PARIS 96}, \cite{sax hjiaj GIENS 97},  
\cite{sax feng IJMCM 98}, \cite{boush chaa IJP 01}, \cite{hjiaj fort IJES 03}, 
\cite{hjiaj feng IJNME 04}). Applications to the contact Mechanics 
\cite{feng hjiaj CM 06}, the Dynamics of granular materials 
(\cite{fort sax CRAS 99},  \cite{fort hjiaj CG 02}, 
\cite{fort mill IJNME 04}\cite{sax fort MOM 04}), the cyclic 
Plasticity of metals \cite{sax hjiaj GIENS 97} and the Plasticity of 
soils (\cite{bersaxce}, \cite{hjiaj fort IJES 03}) illustrate the relevancy of 
this approach.

In all the papers already mentioned about the mechanical applications, 
bipotentials for certain multivalued constitutive laws were constructed. 
Nevertheless, in order to better understand  the 
bipotential approach, one has to solve  the following  problems: 
\begin{enumerate}
\item[1)] (existence) what are the conditions to be satisfied by a multivalued 
 law such that it can be expressed with the help of a bipotential? 
\item[2)] is there a procedure to construct a class of bipotentials for a 
multivalued law? We expect that generically the law does not uniquely determine 
the bipotential. 
\end{enumerate} 

We give a first  mathematical treatment of these  problems and we  prove 
results of existence (theorem \ref{thm1}) and construction (theorem \ref{thm2}) 
 of bipotentials for a class of graphs of multivaluate laws. 

One of the key ideas is constructing the bipotential as an inferior envelope.
 That could be considered as paradoxal because, in general, it is strongly 
unprobable that an inferior envelope, even of convex functions, would be 
convex. 
Nevertheless, we were convinced of the relevancy of this approach by examples 
inspired from mechanics and we wished to understand the reason. 
 That led us to introduce the main tool of convex lagrangian covers 
(Definition \ref{defcover}) satisfying an implicit convexity condition. 

The recipe that we give in this paper applies  only to BB-graphs 
(Definition \ref{dh1}) admitting  at least one convex lagrangian cover by 
maximal cyclically monotone graphs. 
This is an interesting  class of graph of multivalued laws for the 
following two reasons: 
\begin{enumerate}
\item[(a)] it contains the class of graphs of subdifferentials of  convex lsc 
 superpotentials,
\item[(b)] any of the  graphs of non associated laws from the mentioned 
mechanical applications of bipotentials is a BB-graph and it admits a 
physically relevant convex lagrangian cover by cyclically monotone graphs. 
\end{enumerate}    
Relating to point (b), it is important to know that the results from this paper 
don't apply to some  BB-graphs of mechanical interest, 
such as the graph of the  bipotential associated to contact with friction 
\cite{saxfeng}. This is because we use in this paper only convex lagrangian 
covers with {\it maximal} cyclically   monotone graphs, see also 
Remark \ref{ernst1}. 

This paper is only  a first step into the subject of constructions of 
bipotentials.  Our aim is to explain a general method of construction in a 
reasonably simple situation, interesting in itself, leaving aside for the 
moment certain difficulties appearing in  the general method.  Another  article,
 in preparation, is dedicated to the  extension of  the method presented here 
to a more general class of BB-graphs, by relaxing the notion of convex lagrangian 
cover. In this way we shall be able to construct bipotentials even for some of the  BB-graphs described in 
Remark \ref{ernst1}.

\paragraph{Aknowledgements.} The first author  acknowledges partial support 
from the Romanian Ministry of Education and Research, through the grant 
CEEX06-11-12/2006. Part of this work has been done in 2005, when the first 
author has been invited at the Laboratoire de  M\'ecanique de Lille, 
UMR CNRS 8107, Universit\'e des Sciences et Technologies de Lille.  

The authors  thanks E. Ernst  for pointing out an error in a first version 
of the paper,  and for some  examples related to Remarks \ref{ernst1} 
and \ref{ernst2}. We thank also the anonymous referee for comments 
and suggestions leading hopefully to a better paper. 

\section{Notations and Definitions}

$X$ and $Y$ are topological, locally convex, real vector spaces of dual 
variables $x \in X$ and $y \in Y$, with the duality product 
$\langle \cdot , \cdot \rangle : X \times Y \rightarrow \mathbb{R}$. 
We shall suppose that $X, Y$ have topologies compatible with the duality 
product, that is: any  continuous linear functional on $X$ (resp. $Y$) 
has the form $x \mapsto \langle x,y\rangle$, for some $y \in Y$ (resp. 
$y \mapsto \langle x,y\rangle$, for some  $x \in X$). 

For any convex and closed set $A \subset X$, its  indicator function,  
$\displaystyle \chi_{A}$, is defined by 
$$\chi_{A} (x) = \left\{ \begin{array}{ll}
0 & \mbox{ if } x \in A \\ 
+\infty & \mbox{ otherwise } 
\end{array}Ê\right. $$
The indicator  function is convex and lower semi continuous. 
If the set $A$ contains only one element $A=\left\{a\right\}$ then we shall use 
the notation $\displaystyle \chi_{a}$ for the indicator  function of $A$.

We use the notation: $\displaystyle \bar{\mathbb{R}} = \mathbb{R}\cup \left\{ +\infty \right\}$. 

Given a function 
$\displaystyle \phi: X \rightarrow  \bar{\mathbb{R}}$, 
the polar $\phi^{*}: Y \rightarrow \bar{\mathbb{R}}$ 
is defined by: 
$$\phi^{*}(y) = \sup \left\{ \langle y,x\rangle - 
\phi(x) \mid x \in X \right\} \ .$$ 
The polar is always convex and lower semi continuous. 

We denote by $\Gamma(X)$ the class of convex and lower semicontinuous 
functions $\displaystyle \phi: X \rightarrow \bar{\mathbb{R}}$. 
The  class of convex and lower semicontinuous functions 
$\displaystyle \phi: X \rightarrow \mathbb{R}$ is denoted by 
$\displaystyle \Gamma_{0}(X)$.

The subgradient of a function $\displaystyle \phi: X \rightarrow \bar{\mathbb{R}}$ in a point $x \in X$ is the (possibly empty) set: 
$$\partial \phi(x) = \left\{ u \in Y \mid \forall z \in X  \  \langle z-x, u \rangle \leq \phi(z) - \phi(x) \right\} \  .$$ 
In a similar way is defined the subgradient of a function $\psi: Y \rightarrow \bar{\mathbb{R}}$ in a point $y \in Y$, as the set: 
$$\partial \psi(y) = \left\{ v \in X \mid \forall w \in Y  \  \langle v, w-y \rangle \leq \psi(w) - \psi(y) \right\} \ .$$ 

With these notations we have the Fenchel inequality: let 
$\displaystyle \phi: X \rightarrow \bar{\mathbb{R}}$ be a convex lower 
semicontinuous function. Then: 
\begin{enumerate}
\item[(i)] for any $x \in X , y\in Y$ we have $\displaystyle \phi(x) + \phi^{*}(y)  \geq \langle x, y \rangle$; 
\item[(ii)]  for any $(x,y) \in X \times Y$ we have the equivalences: 
$$ y \in \partial \phi(x) \ \Longleftrightarrow \ x \in \partial \phi^{*}(y)  \ \Longleftrightarrow \  \phi(x) + \phi^{*}(y)  = 
\langle x , y \rangle \ . $$
\end{enumerate}

\begin{defi} We model  {\bf the graph of a constitutive law} by a set 
$M \subset X \times Y$. Equivalently, the law is given by the multivalued 
application 
$$\displaystyle  X \ni x  \mapsto m(x) \ = \ \left\{ y \in Y \mid (x,y) \in 
M \right\} \ .$$
The {\bf dual} law is the multivalued application
$$\displaystyle  Y \ni y  \mapsto m^{*}(y) \ = \ \left\{ x \in X \mid (x,y) \in M \right\} \ .$$
The {\bf domain} of the law is the set $\displaystyle dom(M) = \left\{ x \in 
X \mid m(x) \not = \emptyset\right\}$. 
The {\bf image} of the law is the set $\displaystyle im(M) = \left\{ 
y \in Y \mid m^{*}(y) \not = \emptyset\right\}$. 
\label{def1}
\end{defi}

For example, if $\phi: X \rightarrow \mathbb{R}$ is a convex lower semi continuous  function, 
the associated law is the multivalued application 
$\displaystyle \partial \phi$, the subdifferential  of $\phi$,  \cite{moreau} 
Def. 10.1, that is the set of subgradients.The dual law is 
$\partial \phi^{*}$ (the subdifferential of the Legendre-Fenchel dual of $\phi$) 
and the graph of the law is the set 
\begin{equation}
M(\phi) \ = \ \left\{ (x,y) \in X \times Y \ \mid \ \phi(x)+\phi^{*}(y) = \langle x, y \rangle \right\} \  .
\label{mphi}
\end{equation}
For any  convex lower semi continuous function $\phi$ the graph $M(\phi)$ is {\it maximal cyclically 
monotone} (\cite{rocka} Theorem 24.8. or \cite{moreau} Proposition 12.2). Conversely, if 
$M$ is closed and maximal cyclically monotone then there is a convex, 
lower semicontinuous $\phi$ such that 
$M=M(\phi)$. 

\begin{defi} A {\bf bipotential} is a function $b: X \times Y \rightarrow
 \bar{\mathbb{R}}$, with the properties: 
\begin{enumerate}
\item[(a)] $b$ is convex and lower semicontinuos in each argument; 
\item[(b)] for any $x \in X , y\in Y$ we have $\displaystyle b(x,y) \geq \langle x, y \rangle$; 
\item[(c)]  for any $(x,y) \in X \times Y$ we have the equivalences: 
\begin{equation}
y \in \partial b(\cdot , y)(x) \ \Longleftrightarrow \ x \in \partial b(x, \cdot)(y)  \ \Longleftrightarrow \ b(x,y) = 
\langle x , y \rangle \ .
\label{equiva}
\end{equation}
\end{enumerate}
The {\bf graph} of $b$ is 
\begin{equation}
M(b) \ = \ \left\{ (x,y) \in X \times Y \ \mid \ b(x,y) = \langle x, y \rangle \right\} \  .
\label{mb}
\end{equation}
\label{def2}
\end{defi}

{\bf Examples.} {\bf (1.)} (Separable bipotential) To any convex 
lower semicontinuous function  $\phi$  we can associate the 
{\bf separable bipotential} 
$$\displaystyle b(x,y)  = \phi(x) + \phi^{*}(y) .$$ 
The bipotential $b$ and the potential $\phi$ define the same 
law: $\displaystyle M(b) = M(\phi)$. 

{\bf (2.)} (Cauchy bipotential) Let $X=Y$ be a Hilbert space and let the duality 
product be  equal 
to the scalar product. Then we define the {\bf Cauchy bipotential} by the formula 
$$\displaystyle b(x,y) = \| x\| \  \| y\| .$$ 
Let us check the Definition (\ref{def2}) The point (a) is obviously satisfied. 
The point (b) is true by the Cauchy-Schwarz-Bunyakovsky inequality.  
We have equality in the Cauchy-Schwarz-Bunyakovsky inequality 
$b(x,y) = \langle x,y \rangle$ if and only if there is $\lambda > 0$ such 
that $y = \lambda x$ or one of $x$ and $y$ vanishes.  This is exactly the statement from the  point (c), for 
the function $b$ under study. 

\begin{rk}
 The Cauchy bipotential is an ingredient in the construction of many 
bipotentials of mechanical interest,  because the (graph of the) law 
associated to $b$ is the set of pairs of collinear and with same orientation 
vectors. It can not be expressed by a separable potential because $M(b)$ is not 
a cyclically monotone graph. We shall apply the results of this 
paper to the Cauchy bipotential, in order to show that we are able to 
recover the expression of this bipotential from the graph of its associated law.
\label{cauchyrk}
\end{rk}

\section{Existence of a bipotential}

Given a  non empty set $M \subset X \times Y$, Theorem \ref{thm1}
 provides  a necessary and sufficient condition on $M$ for the 
existence of a bipotential $b$ with $M = M(b)$. In 
order to shorten  the notation we shall give a name to this condition: 

\begin{defi}  The non empty set $M \subset X \times Y$ is a {\bf BB-graph } 
(bi-convex, bi-closed) if for all $x \in \ dom(M)$ and for all 
$y \in \ im(M)$ the sets $\displaystyle m(x)$ and 
$\displaystyle m^{*}(y)$ are convex and closed.
\label{dh1}
\end{defi}

The existence problem is easily settled by the following result. 

\begin{thm}
 Given  a non empty set $M \subset X \times Y$, there is a bipotential $b$ 
such that $M=M(b)$ if and only if $M$ is a BB-graph. 
 \label{thm1}
 \end{thm}

\paragraph{Proof.} 
Let $b$ be a bipotential such that $M(b)$ is not void. We first want to prove that for any  $x \in X$ and $y \in Y$ the sets $m(x)$ and $m^{*}(y)$ are convex and closed. 
  
  Indeed, if $m(x)$ or $m(y)$ are empty or they contain only one element then 
there is nothing to prove. Let us suppose, for example, that  $m(x)$ has more 
than one element. From the convexity and lower semi continuity 
hypothesis on $b$ from Definition \ref{def2}, it follows that $m(x)$ is closed and 
convex. Indeed, remark that $m(x)$ is a sub-level set for a convex and lower semi continuous 
mapping: 
$$m(x) = \left\{ y \in Y \mbox{ : } b(x,y) - \langle x,y\rangle \leq 0 \right\} 
\quad , $$
thus a closed and convex set.

  Let us consider now a non empty set $M\subset X \times Y$ such that for any $x \in X$ and $y \in Y$ the sets $m(x)$ and $m^{*}(y)$ are convex and closed. We define then the function $\displaystyle b_{\infty}:  X \times Y \rightarrow \overline{\mathbb{R}}$ by: 
  $$b_{\infty}(x,y) = \left\{ \begin{array}{ll}
\langle x,y \rangle  & \mbox{ if } (x,y) \in M \\ 
+\infty & \mbox{ otherwise } 
\end{array}\right. $$
We have to prove that $\displaystyle b_{\infty}$ is a bipotential and that $\displaystyle M = M(b_{\infty})$. 
This last claim is trivial, so let us check the points from the Definition \ref{def2}. 
For the point (a) notice that for any fixed $x \in X$ the function 
$\displaystyle b_{\infty}(x,\cdot)$ is the sum of a linear continuous 
function with the indicator  function of $m(x)$. By hypothesis the set 
$m(x)$ is closed and convex, therefore its indicator function is convex 
and lower semicontinuous. It follows that  the function  
$\displaystyle b_{\infty}(x,\cdot)$ is convex and lower semi continuous. 
In the same way we prove that for any fixed $y \in Y$  the function  
$\displaystyle b_{\infty}(\cdot, y)$ is convex and lower semi continuous. 
The points (b) and (c) are trivial by the Definition of the function 
$\displaystyle b_{\infty}$. $\quad \blacksquare$

\vspace{.5cm}

\begin{rk}
The uniqueness of $b$  is {\it not true}. For example, in the case of the 
Cauchy bipotential we have two different bipotentials $b$ and $b_{\infty}$ 
with the same graph. Therefore the graph of the law alone is not  sufficient 
to uniquely define the bipotential. 
\label{mostrk}
\end{rk}

\section{Construction of a bipotential}

Theorem \ref{thm1} does not give a satisfying bipotential for a given
 multivalued constitutive law, because the bipotential $b_{\infty}$ is 
definitely not interesting for applications. 

The most important conclusion of preceding section is contained in the Remark 
\ref{mostrk}: in the hypothesis of Theorem \ref{thm1}, the graph of the law 
 is not sufficient to uniquely construct an associated bipotential. This is 
in contrast with the case of a maximal cyclically monotone graph $M$, when 
by Rockafellar theorem (\cite{rocka} Theorem 24.8.) we have a method to 
reconstruct unambigously the associated separable bipotential 
(see point (a) below). 

In our opinion this is the main reason why  the bipotentials are not 
more often used in applications. Without a recipe for constructing 
the bipotential associated with (the experimental data contained in) 
the graph of a non associated mechanical law, there is little chance that 
one may guess a correct expression for this bipotential. 

We are looking for a method of construction of a bipotential  with 
the following properties: 
\begin{enumerate}
\item[(a)] if the graph $M\subset X \times Y$ is maximal cyclically monotone 
then the constructed bipotential is separable (see Example (1.)), 
\item[(b)] the method applied to the graph associated to the Cauchy bipotential 
allows to reconstruct the named bipotential (as mentioned in Remark \ref{cauchyrk}, 
this bipotential appears in many applications), 
\item[(c)] the method should use only hypothesis related to the graph 
$M\subset X \times Y$.
\end{enumerate}

Relating to point (c), we noticed  that in all applications we were able to 
reconstruct the bipotentials by knowing a little more than the graph 
$M\subset X \times Y$, namely a decomposition: 
$$M  \ = \  \bigcup_{\lambda \in \Lambda} M_{\lambda} \quad .$$
We have to mention that in all applications this decomposition 
stems out from physical considerations. 

Thus we were led to the introduction of  convex lagrangian covers.

\begin{defi} Let $M \subset X \times Y$ be a non empty set.  A {\bf convex lagrangian cover }
of  $M$ is a function   
$\displaystyle \lambda \in \Lambda \mapsto \phi_{\lambda}$ from  $\Lambda$ with 
values in the set  $\Gamma(X)$, with the 
properties:
\begin{enumerate}
\item[(a)] The set $\Lambda$ is a non empty compact topological space, 
\item[(b)] Let $f: \Lambda \times X \times Y \rightarrow \bar{\mathbb{R}}$ be the function defined by 
$$f(\lambda, x, y) \ = \ \phi_{\lambda}(x) + \phi^{*}_{\lambda}(y) .$$
Then for any $x \in X$ and for any $y \in Y$ the functions 
$f(\cdot, x, \cdot): \Lambda \times Y \rightarrow \bar{\mathbb{R}}$ and 
$f(\cdot, \cdot , y): \Lambda \times X \rightarrow \bar{\mathbb{R}}$ are  lower 
semi continuous  on the product spaces   $\Lambda \times Y$ and respectively 
$\Lambda \times X$ endowed with the standard topology, 
\item[(c)] We have 
$$M  \ = \  \bigcup_{\lambda \in \Lambda} M(\phi_{\lambda}) \quad  .$$
\end{enumerate}
\label{defcover}
\end{defi}

\section{On the existence and uniqueness of convex lagrangian covers}

Not 
any BB-graph admits a convex lagrangian cover.  There are at least two sources 
of examples of such BB-graphs, described further. For more considerations along this line see the 
last section of the paper. 

\begin{rk} 
Let $M$ be a BB-graph with the property: for any $\phi$, convex, lower semicontinuous 
function defined on $X$, we have 
$\displaystyle M(\phi) \setminus M \not = \emptyset$. 
Then $M$ does not admit any convex lagrangian cover. 

As an example take any 
convex, lower semicontinuous $\displaystyle \phi:X \rightarrow \bar{\mathbb{R}}$ 
and consider $M \subset M(\phi)$, BB-graph, such that  $M \not =  M(\phi)$. 
Then $M$ has the property described previously, therefore it 
 does not admit any convex lagrangian cover. 
\label{ernst1}
\end{rk}

\begin{rk}
If $M$ is a BB-graph and $A$ is any linear, continuous transformation of 
$X \times Y$ into itself, such that $A(X\times \left\{0\right\} ) \subset X \times \left\{ 0 \right\}$ 
and $A(\left\{0\right\} \times Y) \subset \left\{ 0 \right\} \times Y$, 
then $A(M)$ is also a BB-graph. 
However, it may happen that $M$ admits convex lagrangian covers, but not 
$A(M)$. 

Indeed,  we consider $\displaystyle X=Y=\mathbb{R}$ with natural 
duality and a  $\displaystyle \mathcal{C}^{2}$ function 
 $\displaystyle \phi:X \rightarrow \mathbb{R}$ with derivative $\phi'$ strictly 
increasing. Let us define  $M=M(\phi)$ and $A(x,y)=(x,-y)$. 
The set $A(M)$ has a simple description as the graph   
of $-\phi'$. As $\phi'$ is stricly increasing, for any two different 
 $\displaystyle x_{1},x_{2} \in \mathbb{R}$ and 
$\displaystyle y_{i}= - \phi'(x_{i})$ ($i=1,2$), we have 
$$\langle x_{1}-x_{2},  y_{1}-y_{2} \rangle = ( x_{1}-x_{2})( y_{1}-y_{2}) < 0 \quad .$$
This implies that $A(M)$ has the property described in Remark \ref{ernst1}. 
For if there is a convex, lower semicontinuous 
$\displaystyle \psi:X \rightarrow \bar{\mathbb{R}}$ such that 
$M(\psi) \subset A(M)$ then for any two different 
 $\displaystyle x_{1},x_{2} \in \mathbb{R}$ and $\displaystyle y_{i}\in 
\mathbb{R}$, $i=1,2$,  such that 
$\displaystyle (x_{i}, y_{i}) \in M(\psi)$ we would have 
$$\langle x_{1}-x_{2},  y_{1}-y_{2} \rangle = ( x_{1}-x_{2})( y_{1}-y_{2}) \geq  0 \quad ,$$ 
which leads to contradiction. 
\label{ernst2}
\end{rk}

The bipotential 
$\displaystyle b_{\infty}$ from the proof of Theorem \ref{thm1} does not come 
from a convex lagrangian cover.  There exist   BB-graphs  admitting only 
one  convex lagrangian cover (up to reparametrization), as well as BB-graphs 
which have infinitely many lagrangian covers.  

In conclusion, we think  it is a hard and challenging 
mathematical problem to describe all convex lagrangian covers of a BB-graph. 

\section{Implicit convexity and the main result}

The main result of this paper is Theorem \ref{thm2}, which gives a recipe for 
constructing a bipotential not from the graph $M$ of a multivalued law, but 
from a convex lagrangian cover. Therefore the results in this section apply only 
to BB-graphs admitting at least one convex lagrangian cover. 

 In the next section we shall apply this recipe 
for two convex lagrangian covers of $M(b)$, with $b$ equal to the Cauchy bipotential. 

\begin{rk} We give here a justification for the name "convex lagrangian cover". 
Suppose that for any $\lambda \in \Lambda$ the function  
$\displaystyle \phi_{\lambda}$ is smooth. Then it is well known that 
the graph (of the subdifferential of $\displaystyle \phi_{\lambda}$) 
$\displaystyle M(\phi_{\lambda})$ is a lagrangian manifold in the symplectic 
manifold $X \times Y$ with the canonical symplectic form 
$$\omega\left( (x,y), (x',y')\right) = \langle x, y'\rangle - \langle y,x'\rangle \  $$
Therefore the set $M$ is covered by the family of lagrangian manifolds  $\displaystyle M(\phi_{\lambda})$, $\lambda \in \Lambda$. 
\label{langrk}
\end{rk}

With the help of a convex lagrangian cover we shall define a function $b$. 
We intend to prove that (under a certain condition explained further) the 
function $b$ is a bipotential and that $M=M(b)$. 

\begin{defi}
Let $\displaystyle \lambda \mapsto \phi_{\lambda}$ be a convex lagrangian 
cover of the BB-graph $M$. To the cover we associate the function
$b: X \times Y \rightarrow \mathbb{R} \cup \left\{ + \infty \right\}$ 
by the formula 
$$b(x,y) \ = \ \inf \left\{ \phi_{\lambda}(x)+ \phi_{\lambda}^{*}(y) 
\mbox{ : } \lambda \in \Lambda\right\} \ = \  \inf_{\lambda \in \Lambda}  
f(\lambda, x, y)  \quad  . $$
\label{defrecipe}
\end{defi}

We have to check if the function $b$ has the properties (a), (b), (c) 
from the Definition \ref{def2} of a bipotential. 

\begin{prop}
Let $\displaystyle \lambda \mapsto \phi_{\lambda}$ be a convex lagrangian 
cover of the 
BB-graph $M$ and $b$ given by Definition \ref{defrecipe}.  Then: 
\begin{enumerate}
\item[(a)] for all $(x,y) \in M$ we have $b(x,y) \ = \ \langle x ,y \rangle$. 
\item[(b)] for all $(x,y) \in X \times Y$ we have $b(x,y) \geq  \langle x , y \rangle$ .
\end{enumerate}
\label{p2}
\end{prop}

\paragraph{Proof.} 
For all $\lambda \in \Lambda$ and $(x,y) \in X \times Y$ we have the inequality: 
$$\phi_{\lambda}(x)  +  \phi^{*}_{\lambda}(y) \geq \langle x , y \rangle \ .$$
As a consequence of this inequality and Definition \ref{defrecipe} of the function $b$ we obtain the point (b). 

For proving the point (a) it is enough to show that if $(x,y) \in M$ then 
$b(x,y) \leq \langle x , y \rangle$. But this is true. Indeed, if 
$(x,y) \in M$ then there is a $\lambda \in \Lambda$ such that $\displaystyle 
(x,y) \in M(\phi_{\lambda})$ and thus  
$$\phi_{\lambda}(x)  +  \phi^{*}_{\lambda}(y) \ = \ \langle x , y \rangle \quad . $$ 
From the Definition \ref{defrecipe} it follows that for any $\lambda \in \Lambda$ we have 
$$b(x,y) \ \leq \ \phi_{\lambda}(x)  +  \phi^{*}_{\lambda}(y) $$ 
therefore $b(x,y) \leq \langle x , y \rangle$, which finishes the proof . \quad $\blacksquare$ \vspace{.5cm}

\begin{prop}
Let $\displaystyle \lambda \mapsto \phi_{\lambda}$ be a convex lagrangian cover 
of the BB-graph $M$ and $b$ given by Definition \ref{defrecipe}.  
\begin{enumerate}
\item[(a)] Suppose that $x \in X$ is given and that $y \in Y$ has the minimum property 
$$b(x,y)  -  \langle x , y \rangle \ \leq \ b(x,z)  -  \langle x , z \rangle $$ 
for any $z \in Y$. Then $b(x,y) \ = \ \langle x , y \rangle$. 
\item[(b)]  If  $b(x,y) \ = \ \langle x , y \rangle $ then $(x,y) \in M$. 
\end{enumerate}
\label{p3}
\end{prop}

\paragraph{Proof.} 

(a) We start from the Definition of $b$. We have 
$$b(x,y)   \ = \ \inf \left\{ \phi_{\lambda}(x) + \phi^{*}_{\lambda}(z)  \mbox{ : } \lambda \in \Lambda \right\} \ .$$
We  use the compactness of $\Lambda$ (point (a) from  
Definition \ref{defcover}) to obtain a net $\displaystyle (\lambda_{n})_{n}$ in 
$\Lambda$, which converges to $\bar{\lambda} \in \Lambda$, such that 
$b(x,y)$ is the limit of the net 
$\displaystyle  \left( 
\phi_{\lambda_{n}}(x) + \phi^{*}_{\lambda_{n}}(y) \right)_{n}$.

From the lower semicontinuity of the cover (point (b) from Definition \ref{defcover}) we infer that 
$$b(x,y)   \ = \ \phi_{\overline{\lambda}}(x) \  + \  \phi^{*}_{\overline{\lambda}}(y) \ .$$
Remark that the value of the limit $\bar{\lambda}$ of the net 
$\displaystyle (\lambda_{n})_{n}$  depends on $(x,y)$.

The hypothesis from point (a) and the definition of the function $b$ implies that for any $z \in Y$ and any $\lambda \in \Lambda$ we have 
$$\phi_{\overline{\lambda}}(x) \  + \  \phi^{*}_{\overline{\lambda}}(y) - \langle x , y \rangle \ \leq \ 
\phi_{\lambda}(x) \  + \  \phi^{*}_{\lambda}(z) -  \langle x,z\rangle  . $$
In particular, for $\displaystyle \lambda = \overline{\lambda}$ we get that for all $z \in Y$
$$\phi^{*}_{\overline{\lambda}}(y) - \phi^{*}_{\overline{\lambda}}(z) \ \leq \ \langle x, y-z\rangle .$$
This means that $\displaystyle x \in \partial \phi^{*}_{\overline{\lambda}}(y)$, which implies that 
$$b(x,y) \ = \  \phi_{\overline{\lambda}}(x) \  + \  \phi^{*}_{\overline{\lambda}}(y) \ = \ \langle x , y \rangle .$$

For the point (b), suppose that $b(x,y)  =  \langle x , y \rangle$. 
As we remarked before, there is  
$\overline{\lambda} \in \Lambda$ such that 
$$b(x,y) \ = \ \phi_{\overline{\lambda}}(x) \  + \  \phi^{*}_{\overline{\lambda}}(y) \ . $$
Putting all together we see that 
$$\phi_{\overline{\lambda}}(x) \  + \  \phi^{*}_{\overline{\lambda}}(y) \ = \ \langle  x , y \rangle \ ,  $$
therefore $\displaystyle (x,y) \in M(\phi_{\overline{\lambda}}) \subset M$.  \quad $\blacksquare$ \vspace{.5cm} 

We shall give now a sufficient hypothesis for  the separate convexity of $b$. 
This is the last ingredient that we need in order to prove that $b$ is a
 bipotential. 

We shall use  the following  notion of implicit convexity. 

\begin{defi}
Let $\Lambda$ be an arbitrary non empty set and $V$ a real vector space. The 
function $f:\Lambda\times V \rightarrow \bar{\mathbb{R}}$ is 
{\bf implicitly  convex} if for any two elements 
$\displaystyle (\lambda_{1}, z_{1}) , 
(\lambda_{2},  z_{2}) \in \Lambda \times V$ and for any two numbers 
$\alpha, \beta \in [0,1]$ with $\alpha + \beta = 1$ there exists 
$\lambda  \in \Lambda$ such that 
$$f(\lambda, \alpha z_{1} + \beta z_{2}) \ \leq \ \alpha 
f(\lambda_{1}, z_{1}) + \beta f(\lambda_{2}, z_{2}) \quad .$$
\label{defimpl}
\end{defi}

Let us state the last hypothesis from our construction, as a definition. 

\begin{defi}
Let $\displaystyle \lambda \mapsto \phi_{\lambda}$ be a convex lagrangian cover 
of the BB-graph $M$ and $f: \Lambda \times X \times Y \rightarrow \mathbb{R}$ the associated function introduced in Definition \ref{defcover}, that is the function defined by 
$$f(\lambda, z, y) \ = \ \phi_{\lambda}(z) + \phi^{*}_{\lambda}(y) \quad  .$$ 
The cover is bi-implicitly convex (or a  {\bf BIC-cover}) if 
for any $y \in  Y$ and $x\in X$  the functions $f(\cdot, \cdot, y)$ and 
$f(\cdot, x, \cdot)$ are implicitly convex in the sense of Definition 
\ref{defimpl}. 
\label{defbit}
\end{defi}

In the case of $M=M(\phi)$, with $\phi$ convex and lower semi continuous  
(this corresponds to separable bipotentials), the set $\Lambda$ has only one 
element $\Lambda = \left\{ \lambda \right\}$ 
  and we have only one potential $\displaystyle \phi$. The associated 
bipotential from Definition \ref{defrecipe} is obviously 
$$b(x,y) \ = \ \phi(x) + \phi^{*}(y) \ .$$
This is a BIC-cover in a trivial way: the implicit convexity conditions are equivalent with 
the convexity of $\displaystyle \phi$, $\displaystyle \phi^{*}$ respectively. 

Therefore, in the case of separable bipotentials the BIC-cover condition is 
trivially true.

Our recipe  concerning the construction of a bipotential is based on 
the following result. 

\begin{thm} Let $\displaystyle \lambda \mapsto \phi_{\lambda}$ be a BIC-cover of 
the BB-graph $M$ and $b: X \times Y \rightarrow R$ defined by
\begin{equation}
b(x,y) \ = \ \inf \left\{ \phi_{\lambda}(x) + \phi^{*}_{\lambda}(y) \ \mid \  \lambda \in \Lambda \right\} \ . 
\end{equation}
Then $b$ is a bipotential and $M=M(b)$. 
\label{thm2}
\end{thm}

\paragraph{Proof.} 
(Step 1.) We prove first that for any $x \in X$ and for any 
$y \in Y$,  the functions $b(x,\cdot)$ and $b(\cdot, y)$ are convex. 

For fixed $y \in Y$, for any $\displaystyle x_{1},x_{2} \in X$
 and for any $\varepsilon > 0$, there are $\displaystyle \lambda_{1}, \lambda_{2} 
\in \Lambda$ such that ($i=1,2$) 
$$b(x_{i},y) + \varepsilon \geq f(\lambda_{i},x_{i},y) \quad . $$
For the  pairs $\displaystyle (\lambda_{1},x_{1}), (\lambda_{2},x_{2})$ we use 
the implicit convexity of $f(\cdot,\cdot,y)$ to find that there is  
$\lambda \in \Lambda$ such that 
$$f(\lambda,\alpha x_{1}+\beta x_{2},y) \leq \alpha f(\lambda_{1},x_{1},y)+ \beta f(\lambda_{2},x_{2},y) \quad .$$
All in all we have: 
$$ b(\alpha x_{1} + \beta x_{2}) \leq f(\lambda,\alpha x_{1}+\beta x_{2},y) \leq $$
$$\leq  \alpha f(\lambda_{1},x_{1},y) + \beta f(\lambda_{2},x_{2},y) \leq \alpha b(x_{1}, y) + \beta b(x_{2},y) + \varepsilon \quad . $$

As $\varepsilon >0$ is an arbitrary chosen positive number, the convexity of the function $b(\cdot, y)$ is proven. The proof for the convexity of $b(x, \cdot)$ is similar. 

(Step 2.) We shall prove now that for any $x \in X$ and for any $y \in Y$,  
the functions $b(\cdot, x)$ and $b(\cdot, y)$ are lower semicontinuous. 
Consider a  net  
$\displaystyle (x_{n})_{n}\in X$ which converges to $x$. We use the same reasoning as in the proof of Proposition \ref{p3} (a) to deduce that for  each $n \in \mathbb{N}$ there exists a $\displaystyle \lambda_{n} \in \Lambda$ such that 
$$b(x_{n}, y) \ = \ \phi_{\lambda_{n}}(x_{n}) + \phi_{\lambda_{n}}(y) = f(\lambda_{n}, x_{n}, y) .$$
$\Lambda$ is compact, therefore up to the choice of a subnet, there exists a $\lambda \in \Lambda$ such that $(\lambda_{n})_{n}$ converges to $\lambda$. We use now the lower semicontinuity of $f(\cdot, \cdot, y)$ in order to get that 
$$b(x,y) \leq f(\lambda, x, y) \leq \liminf_{n \rightarrow \infty} f(\lambda_{n}, x_{n}, y) ,$$
therefore the lower semicontinuity of $b(x, \cdot)$ is proven. 
For the function  $b(\cdot, y)$ the proof is similar.  

(Step 3.) $M=M(b)$. Indeed, this is true, by Propositions \ref{p2} (a) and \ref{p3} (b). 

(Step 4.) By Proposition \ref{p2} (b) we have that for any $(x,y) \in X \times Y$ the inequality $b(x,y) \geq \langle x, y \rangle$ is true. Therefore the conditions (a), (b),  from the Definition \ref{def2} of a bipotential, 
are verified.  

(Step 5.) The only thing left to prove is the string of equivalences from Definition \ref{def2} (c). Using the knowledge that $b$ is separately convex and lower semicontinuous, we remark that in fact we only have to prove two implications. 

The first is: for any  $x \in X$ suppose that  $y \in Y$ has the minimum property 
$$b(x,y)  -  \langle x , y \rangle \ \leq \ b(x,z)  -  \langle x , z \rangle $$ 
for any $z \in Y$. Then $b(x,y) \ = \ \langle x , y \rangle$.  

The second implication is similar, only that we start with an arbitrary  $y \in Y$ and with $x \in X$ satisfying the minimum property 
$$b(x,y)  -  \langle x , y \rangle \ \leq \ b(z,y)  -  \langle z , y \rangle $$ 
for any $z \in X$. Then $b(x,y) \ = \ \langle x , y \rangle$. 

The first implication is just Proposition \ref{p3} (a). The second implication 
has a similar proof.   \quad $\blacksquare$ 

The next proposition makes easier to check if a convex 
lagrangian cover satisfies the BIC condition. 

\begin{prop}
Let $\displaystyle \lambda \mapsto \phi_{\lambda}$ be a BIC-cover of 
the BB-graph $M$. Consider 
 any $\alpha,\beta \in [0,1]$, $\alpha+\beta=1$,  any 
$y \in  im(M)$,  any $\displaystyle \lambda_{1}, \lambda_{2} \in \Lambda$ 
and any $\displaystyle x_{1} \in \partial \phi^{*}_{\lambda_{1}}(y)$, 
 $\displaystyle x_{2} \in \partial \phi^{*}_{\lambda_{2}}(y)$. 
According to the BIC condition 
there  exists 
$\displaystyle \lambda \in \Lambda$ such that 
\begin{equation}
f(\lambda, \alpha x_{1} + \beta x_{2}, y) \leq \alpha f(\lambda_{1},x_{1},y) + 
\beta f(\lambda_{2},x_{2},y) \quad .
\label{needd}
\end{equation}
Then $\lambda$ has the property:

$$\alpha  x_{1}  +  \beta x_{2} \ \in  \ \partial \phi_{\lambda}^{*}(y) \ .$$

\label{p1}
\end{prop}

\paragraph{Proof.}

Inequality (\ref{needd}) expresses as: 
\begin{equation}
\phi_{\lambda}(\alpha x_{1} + \beta x_{2}) + \phi_{\lambda}^{*}(y) \leq 
\alpha \phi_{\lambda_{1}}(x_{1}) + \beta \phi_{\lambda_{2}}(x_{2}) + 
\alpha \phi_{\lambda_{1}}^{*}(y) + \beta \phi_{\lambda_{2}}^{*}(y) \quad . 
\label{need1}
\end{equation}
We have also ($i=1,2$)  
$$\phi_{\lambda_{i}}(x_{i}) + \phi_{\lambda_{i}}^{*}(y) = \langle x_{i}, y\rangle \quad .$$
We use this in the inequality (\ref{need1}) to get 
$$\phi_{\lambda}(\alpha x_{1} + \beta x_{2}) + \phi_{\lambda}^{*}(y) \leq 
\langle \alpha x_{1} + \beta x_{2}, y\rangle \quad , $$
which shows that $\displaystyle \alpha x_{1}+\beta x_{2} \in \partial \phi_{\lambda}^{*}(y)$. 
Therefore the $\lambda \in \Lambda$ given by the implicit convexity inequality satisfies 
the conclusion of the proposition.  \quad $\blacksquare$

\begin{rk}
Enforcing the satisfaction of the implicit convexity 
inequality for all  values of $\lambda$ which satisfy the conclusion of 
Proposition \ref{p1}  would be too strong. This remark is 
supported by the second example in section \ref{secrecon}, involving a family 
of non differentiable potentials for which there is no uniqueness for $\lambda$. 
\label{remchoose}
\end{rk}

\section{Reconstruction of the Cauchy bipotential}
\label{secrecon}

In this section we shall reconstruct the Cauchy bipotential from two different 
convex lagrangian covers. As explained in Remark \ref{cauchyrk}, it is important 
for applications that we are able to reconstruct the expression of the Cauchy 
bipotential from 
the graph of its associated law. 

We shall take $X=Y=\mathbb{R}^{n}$ and the duality product is the usual scalar  
product in $\mathbb{R}$. The Cauchy bipotential is 
$\displaystyle \overline{b}(x,y) = \| x\| \|y\|$.  By Cauchy-Schwarz-Bunyakovsky 
inequality the set $M = M(\overline{b})$ is 
$$M = \left\{ (x,y) \in \mathbb{R}^{n} \times  \mathbb{R}^{n} \mbox{ : } \exists \lambda > 0 \ , x=\lambda y\right\} \cup \left( \left\{ 0 \right\} \times \mathbb{R}^{n}\right) \times \left(  \mathbb{R}^{n} \times \left\{ 0 \right\} \right) .$$

Let us consider the topological compact set $\Lambda = [0,\infty]$ 
(with usual topology) and the function 
$\displaystyle \lambda \in \Lambda \mapsto \phi_{\lambda}$ defined as:
\begin{enumerate}
\item[-] if $\lambda \in [0,\infty)$ then 
$\displaystyle \phi_{\lambda}(x) = \frac{\lambda}{2} \| x\|^{2}$, 
\item[-] if $\lambda = \infty$ then 
$$ \phi_{\infty}(x) = \chi_{0} (x) = \left\{ \begin{array}{ll}
0 & \mbox{ if } x = 0 \\ 
+\infty & \mbox{ otherwise } 
\end{array} \right. $$
\end{enumerate}
A straightforward computation shows that the associated function $f$ has the 
expression: 
\begin{equation}
f(\lambda, x, y)  = \left\{ \begin{array}{ll}
\frac{\lambda}{2} \| x\|^{2} + \frac{1}{2\lambda} \|y\|^{2} & \mbox{ if } \lambda \in (0,\infty) \\ 
\chi_{0}(y) & \mbox{ if } \lambda = 0 \\
\chi_{0}(x) & \mbox{ if } \lambda = \infty
\end{array} \right. 
\label{ef1}
\end{equation}

It is easy to check that we have here a convex lagrangian cover of the set $M$. 
We shall prove now that we have a BIC-cover, according to Definition \ref{defbit}. 

The cases $\lambda = 0$ and $\lambda = \infty$ will be treated separately. 

Consider $y \in  im(M) = \mathbb{R}^{n}$, $\displaystyle x_{1}, x_{2} \in X$,  
$\alpha,\beta \in [0,1]$, $\alpha+\beta=1$, and 
$\displaystyle \lambda_{1}, \lambda_{2} \in (0,\infty)$. 
We have to find $\lambda \in \Lambda$ such that 
\begin{equation}
f(\lambda,\alpha x_{1} + \beta x_{2},y) \leq \alpha f(\lambda_{1},x_{1},y) + \beta f(\lambda_{2},x_{2},y) \quad . 
\label{need2}
\end{equation}

We use  Proposition \ref{p1},  for $i=1,2$ and  
$\displaystyle x_{i} \in \partial 
\phi_{\lambda_{i}}^{*}$ in order to find   the value of $\lambda$. 
Computation shows that there is only one such   $\lambda \in \Lambda$, given by 
 \begin{equation}
 \frac{1}{\lambda} \ = \ \frac{\alpha}{\lambda_{1}} + \frac{\beta}{\lambda_{2}} \quad . 
 \label{eqex1}
 \end{equation}
 As this value depends only on $\displaystyle \lambda_{1}, \lambda_{2}$, we shall try to see if 
 this $\lambda$ is good for any choice of $\displaystyle x_{1}, x_{2}$. 
 
 This is indeed the case:   with $\lambda$ given by (\ref{eqex1}) 
the relation (\ref{need2}) (multiplied by $2$)  becomes:  
 \begin{equation}
 \lambda \| \alpha x_{1} + \beta x_{2}\|^{2} \leq \alpha \lambda_{1} \| x_{1}\|^{2} + \beta  \lambda_{2} \| x_{2}\|^{2} \quad . 
\label{need3}
\end{equation}
 Remark that (\ref{eqex1}) can be written as: 
 $$\frac{\alpha \lambda}{\lambda_{1}} + \frac{\beta \lambda}{\lambda_{2}} \ = \ 1 \ . $$
 Write then the fact that the square of the norm is convex, for the convex combination 
 of $\displaystyle \lambda_{1} x_{1} , \lambda_{2} x_{2}$, with the coefficients $\displaystyle \frac{\alpha \lambda}{\lambda_{1}},  \frac{\beta \lambda}{\lambda_{2}} $. We get, after easy simplifications, the 
 inequality (\ref{need3}).

 If $\displaystyle \lambda_{1} = 0$,  
$\displaystyle \lambda_{2} \in (0,\infty)$ then $y$ has to be equal to $0$ and 
$\displaystyle x_{1}$ is arbitrary, $\displaystyle x_{2} = 0$ and 
$\lambda = 0$. The inequality  (\ref{need2}) is then trivial. 
 
 All other exceptional cases lead to trivial inequalities. 
 
 Remark that for any $\lambda \in \Lambda$ and any $x,y \in \mathbb{R}^{n}$ we have 
 $$f(\lambda, x, y) = f(\frac{1}{\lambda}, y, x)$$
 with the conventions $1/0 = \infty$, $1/\infty = 0$. This symmetry and previous proof imply that we have a BIC-cover. 
 
 We compute now the function $b$ from Definition \ref{defrecipe}. We know from Theorem \ref{thm2} that $b$ is a bipotential for the set $M$. 
 
 We have:
 $$b(x,y) = \inf \left\{ f(\lambda, x,y) \mbox{ : } \lambda \in [0,\infty] \right\} .$$
 From the relation (\ref{ef1}) we see that actually 
 $$b(x,y) = \inf \left\{ \frac{\lambda}{2} \| x\|^{2} + \frac{1}{2\lambda} \|y\|^{2}  \mbox{ : } \lambda \in (0,\infty) \right\} .$$
 By the arithmetic-geometric mean inequality we obtain that $b(x,y) = \| x\| \|y\|$, that is the Cauchy bipotential. 
 
 Here is a second example, supporting  the  Remark  \ref{remchoose}.  We shall 
reconstruct the Cauchy bipotential starting 
from a family of non differentiable convex potentials. 
 
 Let $\lambda \geq 0$ be non negative and the closed ball of center $0$ and 
radius $\lambda$ be defined by 
 $$B(\lambda) \ = \ \left\{ y \in Y \mbox{ : } \| y \| \leq \lambda \right\} \ \ \ .$$
 Defining $B(+\infty)$ as the whole space $Y$, one can suppose that 
$\lambda$ belongs to the compact set $\Lambda = [0,+\infty]$. 
 
 For $\lambda \in [0,+\infty)$ we define the set: 
 $$M_{\lambda} \ = \ \left\{ (0,y) \in X\times Y \mbox{ : } \| y \| < \lambda \right\} \cup \left\{ (x,y) \in X\times Y 
 \mbox{ : } \|y\| = \lambda \mbox{ and } \exists \eta \geq 0 \ \ x = \eta y \right\} \ .$$
  One can recognize $\displaystyle M_{\lambda}$  as the graph of the yielding law of a plastic material with a yielding  threshold equal to $\lambda$. 
 For $\lambda = + \infty$ we set 
$\displaystyle M_{+\infty} \ = \  \left\{ 0\right\} \times Y$.
 
 It can be easily verified that the family 
$\displaystyle (M_{\lambda})_{\lambda \in \Lambda}$ of maximal cyclically 
monotone  graphs provides us   a convex lagrangian cover of the set: 
 $$M \ = \ \left\{ (x,y) \in X\times Y \mbox{ : } \exists \alpha, \beta \geq 0 \ \  \alpha x =  \beta y \right\} \ .$$
 The corresponding convex lagrangian cover  is given by: 
 \begin{enumerate}
 \item[-] for $\lambda \in [0,+\infty)$, $\displaystyle \phi_{\lambda}(x) \ = \ \lambda \| x \|$, $\displaystyle 
 \phi_{\lambda}^{*}(y) \ = \ \chi_{B(\lambda)}(y)$, 
 \item[-] $\displaystyle \phi_{+\infty} (x) \ = \ \chi_{0}(x)$, $\displaystyle \phi_{+\infty}^{*} (y) \ = 0$. 
 \end{enumerate}
 The associated function $f$ has the expression: 
\begin{equation}
f(\lambda, x, y)  = \left\{ \begin{array}{ll}
\lambda \| x\|+ \chi_{B(\lambda)}(y) & \mbox{ if } \lambda \in (0,\infty) \\ 
\chi_{0}(y) & \mbox{ if } \lambda = 0 \\
\chi_{0}(x) & \mbox{ if } \lambda = +\infty
\end{array} \right. 
\label{ef11}
\end{equation}
 
 All hypothesis excepting the BIC-cover condition  are obviously satisfied. 
We check this condition  further. 
 Let $\displaystyle \lambda_{1} < \lambda_{2}$, both in $[0,+\infty)$. We want 
first  to determine the values of $\lambda$ fulfilling the conclusion 
 of Proposition \ref{p1}. Let us recall that: 
 \begin{enumerate}
 \item[-] if $\|y\| < \lambda$ then $\displaystyle \partial \phi^{*}_{\lambda}(y) \ = \ \left\{ 0 \right\}$, 
  \item[-] if $\|y\| = \lambda$ then $\displaystyle x \in \partial \phi^{*}_{\lambda}(y)$ is equivalent to: $\exists \eta \geq 0$ such that $x = \eta y$,
   \item[-] if $\|y\| > \lambda$ then $\displaystyle \partial \phi^{*}_{\lambda}(y) \ = \ \emptyset$.
   \end{enumerate}
   
  Then the following events have to be considered: 
  \begin{enumerate}
  \item[(1)] if $\displaystyle \|y \| < \lambda_{1} < \lambda_{2}$  then $\displaystyle x_{1} \in \partial \phi_{\lambda_{1}}(y)$ and $\displaystyle x_{2} \in \partial \phi_{\lambda_{2}}(y)$ imply $\displaystyle 
  x_{1}=x_{2}=0$, 
  \item[(2)] if $\displaystyle \|y \| = \lambda_{1} < \lambda_{2}$  then $\displaystyle x_{1} \in \partial \phi_{\lambda_{1}}(y)$ and $\displaystyle x_{2} \in \partial \phi_{\lambda_{2}}(y)$ imply: $\exists \eta \geq 0$ such that  $\displaystyle 
  x_{1}= \eta y $ and $\displaystyle x_{2}=0$. Thus 
  $$\alpha x_{1} + \beta x_{2} \ = \ \alpha \eta x_{1} \in \partial \phi^{*}_{\lambda}(y)$$ 
  occurs for any $\lambda \geq \|y\|$ when $\displaystyle x_{1}=0$ and $\lambda = \|y\|$ otherwise. 
  \item[(3)] If $\displaystyle \lambda_{1} < \|y\|$ then there is no $\displaystyle x_{1}$ such that $\displaystyle x_{1} \in \partial \phi^{*}_{\lambda_{1}}(y)$. Likewise, if $\displaystyle \lambda_{2} < \|y\|$ then there is no $\displaystyle x_{2}$ such that $\displaystyle x_{2} \in \partial \phi^{*}_{\lambda_{2}}(y)$.
 \end{enumerate}

Consider $y \in  im(M) = \mathbb{R}^{n}$, $\displaystyle x_{1}, x_{2} \in X$,  $\alpha,\beta \in [0,1]$, $\alpha+\beta=1$, and $\displaystyle \lambda_{1}, \lambda_{2} \in [0,\infty)$. For the verification of the implicit convexity inequality (\ref{need2}), we 
need only to consider the case $\displaystyle \| y \| \leq \min 
\left\{ \lambda_{1}, \lambda_{2} \right\}$ and we shall choose  $\lambda= \min 
\left\{ \lambda_{1}, \lambda_{2} \right\} \geq \|y\|$. 
The relation (\ref{need2}) becomes 
$$ \min\left\{ \lambda_{1}, \lambda_{2} \right\} \| \alpha x_{1}+\beta x_{2} \| \leq \alpha \lambda_{1} \| x_{1}\| + \beta \lambda_{2} \|x_{2}\|  \quad , $$
which is true by the convexity of the norm.  All the other cases turn out to be trivial. 

The other half of the BIC-cover condition has a similar proof 
(remark though that the associated function $f$ is not symmetric, 
as in the previous case). 
 
 By virtue of Theorem \ref{thm2},  the function given by Definition \ref{defrecipe}, namely 
 $$b(x,y) \ = \ \inf \left\{ \phi_{\lambda}(x)+ \phi_{\lambda}^{*}(y) \mbox{ : } \lambda \in [0,+\infty] \right\} \  , $$
 is a bipotential. Computation shows that $b$ is the Cauchy bipotential. Indeed: 
 $$b(x,y)  \ = \ \inf \left\{ \lambda \|x\| + \chi_{B(\lambda)}(y) \mbox{ : } \lambda \in [0,\infty) \right\} \ = \ $$
 $$= \ \inf \left\{ \lambda \| x \| \mbox{ : } \lambda \geq \| y \| \right\} \ = \ \|y\| \|x\| \ .$$

\section{Conclusion and perspectives}

Given (the graph of) a multivalued constitutive law $M$,  there is a  bipotential $b$ such that 
$M=M(b)$  if and only if $M$ is a  BB-graph (Definition \ref{dh1} and Theorem \ref{thm1}).   If  
the BB-graph $M$ admits a convex lagrangian cover (Definition \ref{defcover}) which is bi-implicitly 
convex (Definition \ref{defbit}) then we are able to construct an associated  
 bipotential (Theorem \ref{thm2}). 
 
 Remarks \ref{ernst1} and \ref{ernst2}     show that not any BB-graph admits a convex lagrangian cover. 
 We would like to elaborate on the obstructions to the existence of such covers. We start with the 
 example from the  Remark  \ref{ernst2}, due to E.  Ernst.  
 
 From a mechanical point of view, multivalued laws  $M$ with the property that for any two different 
 pairs $\displaystyle (x_{1},  y_{1}), (x_{2}, y_{2}) \in M$ we have 
 $$\langle x_{1} - x_{2}, y_{1} - y_{2}\rangle < 0 $$
 are not  very interesting. Indeed, suppose that the evolution of a mechanical system is described by a 
 sequence of states $\displaystyle (x_{n}, y_{n}) \in M$. Then, as the system passes from one state to another, the work done is always negative. Much more interesting seem to be multivalued laws with the 
 property that for any $\displaystyle (x,y) \in M$ there is at least a different pair $\displaystyle (x',y')\in M$ 
 such that 
 $$\langle x-x', y-y'\rangle \geq 0 \quad .$$
The BB-graphs admitting a convex lagrangian cover have this property. 
 
 There is another aspect, concerning the linear transformation $A$ from the 
Remark \ref{ernst2}. In the example given the transformation $A(x,y) = (x,-y)$ 
is not symplectic, but still it transforms lagrangian sets into lagrangian sets.
  In general, if the dimension of $X$ is strictly greater than one then we 
 can find linear endomorphisms of $X\times Y$ transforming  lagrangian subsets 
of $X\times Y$ into  sets which are not lagrangian, thus destroying lagrangian 
covers.  Moreover, we can find linear symplectic transformations  which 
transforms a convex lagrangian cover into a lagrangian cover which is no longer 
convex. For example,  take $\displaystyle X=Y=\mathbb{R}$, $A(x,y)=(x,y-x)$ and 
the BB-graph $\displaystyle M= \mathbb{R}\times \left\{0\right\}$. 
Then $det(A)=1$, therefore $A$ is symplectic, and 
 $ \displaystyle A(M) = \left\{ (x,-x) \mbox{ : } x \in \mathbb{R} \right\}$. 
The set $A(M)$ is a BB-graph and  a lagrangian set, but it does not admit a 
convex lagrangian cover. The reason for this phenomenon is that convexity is 
not a symplectic invariant. Nevertheless, there are famous theorems in 
Hamiltonian Dynamics which have 
 a convexity assumption in the hypothesis, like the theorem of  Rabinowitz 
stating that the Reeb vector field on the  boundary of a convex domain which is 
bounded  has at least a closed orbit 
 (equivalently, a convex and coercive hamiltonian on $\mathbb{R}^{2n}$ admits a 
closed orbit on every level set). We can easily destroy the convexity 
assumption of this theorem but not the conclusion, by applying a nonlinear 
symplectomorphism.   
 
 For the notion of convex lagrangian cover we had the following source of inspiration. If $M$ is  a symplectic manifold and with 
 convexity assumptions left aside, lagrangian covers as described in this paper resemble to  (real) symplectic polarizations, which are a basic tool in some problems of symplectic geometry. 
 
 Much more interesting are cases relating with Remark \ref{ernst1}. We may consider BB-graphs $M$ 
 not admitting convex lagrangian covers, but with the property that there is a family of convex, lower 
 semicontinuous functions $\displaystyle \phi_{\lambda}$, $\lambda \in \Lambda$ such that 
 $$M \subset \bigcup_{\lambda \in \Lambda} M(\phi_{\lambda}) $$
 with strict inclusion.
 This is the case, for example, of the bipotential which appears in \cite{saxfeng}, related to contact with friction. In a future  paper we shall  extend this method of convex lagrangian cover to 
lagrangian covers by graphs which are cyclically monotone but not 
necessarily maximal cyclically monotone.

\vspace{\baselineskip}

\end{document}